\documentclass[12pt]{article}
\usepackage{amssymb}
\usepackage{amsmath}

\begin{document}

\begin{center}
\textbf{\large Sets with distinct sums of pairs, long
arithmetic progressions, and continuous mappings}
\end{center}

\begin{center}
\textsc{V. Lebedev}
\end{center}

\begin{quotation}
{\small \textbf{Abstract:} We show that if $\varphi \colon
\mathbb R\rightarrow\mathbb R$ is a continuous mapping and the
set of nonlinearity of $\varphi$ has nonzero Lebesgue measure,
then $\varphi$ maps bijectively a certain set that contains
arbitrarily long arithmetic progressions onto a certain set
with distinct sums of pairs.

\textbf{Key words:} thin sets, sets with distinct sums of
pairs, sets of type $\Lambda(4)$, long arithmetic progressions.

2010 Mathematics Subject Classification: Primary 43A46,
Secondary 42A75.}
\end{quotation}

\quad

\begin{center}
\textbf{1. Introduction}
\end{center}

We say that a set $E$ in the real line $\mathbb R$ is a set
with distinct sums of pairs if a relation
$\lambda_1+\lambda_2=\lambda_3+\lambda_4$ with $\lambda_j\in
E$, $1\leq j\leq 4$, holds only in the trivial case when
$\lambda_1=\lambda_3$ and $\lambda_2=\lambda_4$ or when
$\lambda_1=\lambda_4$ and $\lambda_2=\lambda_3$.

By an arithmetic progression of length $N$ we mean a set
$F\subseteq\mathbb R$ of the form $F=\{x+ky, \,k=1, 2, \dots,
N\}$, where $x, y\in\mathbb R$ and $y\neq 0$.

It is natural to regard sets with distinct sums of pairs as
thin sets. One of the properties of sets with distinct sums of
pairs is as follows. Let $B^p(\mathbb R), 1\leq p\leq\infty,$
denote the Besicovitch spaces of almost periodic functions (the
definition of the spaces $B^p(\mathbb R)$ and basic facts about
them can be found in [1] and [10]). Recall that the norm
$\|\cdot\|_{B^p(\mathbb R)}$ on $B^p(\mathbb R)$ is defined by
$$
\|f\|_{B^p(\mathbb R)}=
\limsup_{T\rightarrow+\infty}\bigg(\frac{1}{2T}\int_{-T}^T |f(t)|^p dt\bigg)^{1/p}.
$$
It is well known that for $1\leq p_1\leq p_2\leq\infty$ we have
$B^{p_2}(\mathbb R)\subseteq B^{p_1}(\mathbb R)$ with a natural
relation for the norms: $\|\cdot\|_{B^{p_1}(\mathbb
R)}\leq\|\cdot\|_{B^{p_2}(\mathbb R)}$. In particular,
$\|\cdot\|_{B^2(\mathbb R)}\leq\|\cdot\|_{B^4(\mathbb R)}$.
Assume now that $E$ is a set with distinct sums of pairs, and
$f$ is an $E$-polynomial, i.e., a function of the form
$$
f(t)=\sum_{\lambda\in E}c(\lambda)e^{i\lambda t},
$$
where only finitely many coefficients $c(\lambda)$ are nonzero.
Then we have
$$
\|f\|_{B^4(\mathbb R)}\leq c\|f\|_{B^2(\mathbb R)},
\eqno(1)
$$
where $c>0$ does not depend on $f$. This estimate is nearly
obvious, it suffices to note that
$$
|f(t)|^2=\sum_{\lambda_1, \, \lambda_2\in E}c(\lambda_1)\overline{c(\lambda_2)}e^{i(\lambda_1-\lambda_2)t}
$$
and
$$
|f(t)|^4=\sum_{\lambda_1, \,\lambda_2, \,\lambda_3 \,\lambda_4\in E}
c(\lambda_1)c(\lambda_2)\overline{c(\lambda_3)}\,\overline{c(\lambda_4)}
e^{i(\lambda_1+\lambda_2-\lambda_3-\lambda_4)t}.
$$
Since
$$
\lim_{T\rightarrow+\infty}\frac{1}{2T}\int_{-T}^{T}e^{iat} dt=
\left\{\begin{array}{rl} 0 & \text{if}\quad a\neq 0,\\
1 & \text{if}\quad a=0,\end{array} \right.
$$
it follows that
$$
\|f\|_{B^2(\mathbb R)}=\bigg(\sum_\lambda |c(\lambda)|^2\bigg)^{1/2}
$$
and
$$
\|f\|_{B^4(\mathbb R)}\leq 2^{1/4}\bigg(\sum_\lambda |c(\lambda)|^2\bigg)^{1/2},
$$
which yields (1) with $c=2^{1/4}$.

On the other hand, we regard sets that contain arbitrarily long
arithmetic progressions as very massive ones. Let
$$
\gamma_p(N)=\bigg\|\sum_{k=1}^N e^{ikt}\bigg\|_{L^p(\mathbb T)},
$$
where $\mathbb T=\mathbb R/(2\pi\mathbb Z)$ is the circle
($\mathbb Z$ is the additive group of integers). Assume that
$\{a+kd, k=1, 2, \dots, N\}$ is a progression of length $N$
contained in $E$. Consider the polynomial $f_{N, a,
d}(t)=\sum_{k=1}^N e^{i(a+kd)t}$. Note that $\|f_{N, a,
d}\|_{B^p(\mathbb R)}=\gamma_p(N)$. At the same time
$\gamma_4(N)\simeq N^{3/4}$, $\gamma_2(N)\simeq N^{1/2}$
(clearly, $\gamma_p(N)$  behaves as the $L^p(\mathbb T)$-norm
of the Dirichlet kernel $D_N$ for large $N$). So, if $E$
contains arbitrarily long arithmetic progressions then, in
general, estimate (1) for $E$ -polynomials $f$ does not hold.

We note that, since for a $2\pi$~-periodic function its
$B^p(\mathbb R)$~-norm coincides with the $L^p(\mathbb
T)$~-norm, it follows that if a set $E$ with distinct sums of
pairs is in $\mathbb Z$ then (1) has the form
$\|f\|_{L^4(\mathbb T)}\leq c \|f\|_{L^2(\mathbb T)}$ (for any
$E$-polynomial $f$), so, $E$ is a set of type $\Lambda(4)$ (see
[11] for basic results on sets of type~$\Lambda(p)$).

In this paper we consider continuous mappings $\varphi \colon
\mathbb R\rightarrow\mathbb R$ and show that, with a possible
exception for the case when $\varphi$ is of a very special form
resembling that of a piecewise linear mapping, every continuous
mapping $\varphi$ is singular in the sense that it maps
bijectively a certain set that contains arbitrarily long
arithmetic progressions (a massive set) onto a certain set with
distinct sums of pairs (a thin set).

Our interest to the question of how thin the images of
arithmetic progressions under continuous mappings can be is
inspired by the proof of Beurling and Helson of their theorem
[2] on the endomorphisms of the algebra of Fourier transforms
of measures, and the work of Graham [3] on the mappings that
preserve Sidon sets. We discuss these works in remarks at the
end of this paper.

\quad

\begin{center}
\textbf{2. Statement of the result}
\end{center}

We say that $t\in\mathbb R$ is a point of nonlinearity of a
mapping $\varphi\colon  \mathbb R\rightarrow\mathbb R$ if $t$
has no neighborhood in which $\varphi$ coincides with a linear
function.\footnote{For real functions of one variable, we use
the terms ``affine mapping'' and ``linear function'' as
synonyms.} The set of all such points is called the set of
nonlinearity of $\varphi$ and is denoted by $E(\varphi)$.
Clearly, the set $E(\varphi)$ is closed.

The result of this paper is the following theorem.

\quad

\textsc{Theorem.} \emph{Let $\varphi$ be a continuous
self-mapping of $\mathbb R$. Suppose that $E(\varphi)$ has
nonzero Lebesgue measure. Then there exist a set $A$ that
contains arbitrarily long arithmetic progressions and a set $B$
with distinct sums of pairs such that $\varphi$ maps
bijectively $A$ onto $B$.}

\quad

Note that the case when $E(\varphi)$ is finite is trivial.
Assume that $\varphi$ is piecewise linear and maps bijectively
$A$ onto $B$. Then, if $A$ contains arbitrarily long arithmetic
progressions, $B$ has the same property as well. The general
case of mappings whose sets of nonlinearity have measure zero
seems to be a difficult one. In particular the author does not
know whether the classical Cantor staircase function is
singular in the above arithmetical sense.

\quad

\begin{center}
\textbf{3. Statement of the main Lemma and deduction of the
Theorem}
\end{center}

Let $V$ be the family of the following seven vectors in
$\mathbb R^4$:
$$
v^1=(1, 1, -1, -1), \quad v^2=(1, 1, -2, 0), \quad v^3=(1, -1, 0, 0),
\quad v^4=(1, 1, 0, 0),
$$
$$
v^5=(2, 0, 0, 0), \quad v^6=(1, 1, -1, 0),
\quad v^7=(1, 0, 0, 0).
\eqno(2)
$$

\noindent By $\mathbb Z_0^4$ we denote the set of all vectors
in $\mathbb R^4$ with pairwise different integer coordinates,
and by $[1, N]^4$ the cube in $\mathbb R^4$ formed by all
vectors $x=(x_1, x_2, x_3, x_4)\in\mathbb R^4$ satisfying
$1\leq x_j\leq N$, $1\leq j\leq 4$.

The following lemma is the key assertion to the proof of the
theorem.

\quad

\textsc{Lemma 1} (the main Lemma). \emph{Let $\varphi \colon
\mathbb R\rightarrow\mathbb R$ be a continuous mapping. Let
$S\subset\mathbb R$ be a finite set and $N\geq 4$ be an
integer. Suppose that for every $x, y\in\mathbb R$ there exist
a vector $v=(v_1, v_2, v_3, v_4)\in V$ and a vector $k=(k_1,
k_2, k_3, k_4)\in [1, N]^4\cap \mathbb Z_0^4$ such that
$$
\sum_{j=1}^4v_j\varphi(x+k_jy)\in S.
$$
Then $E(\varphi)$ has Lebesgue measure zero.}

\quad

We postpone the proof of this lemma to the next
section; here we show how to derive the theorem from it.

Note that if $E(\varphi)$ has nonzero measure and $N\geq 4$ is
an integer, then, setting $S=\{0\}$ and applying Lemma~1, we
can find $x, y\in\mathbb R$, such that
$$
\sum_{j=1}^4v_j\varphi(x+k_jy)\neq 0
$$
for all vectors $v=(v_1, v_2, v_3, v_4)\in V$ and $k=(k_1, k_2,
k_3, k_4)\in [1, N]^4\cap \mathbb Z_0^4$. Since the family $V$
contains the vectors $v^1, v^2$ and $v^3$ (see (2)), it follows
that the numbers $\varphi(x+ky), \,k=1, 2, \ldots, N,$ are
pairwise distinct and form a set with distinct sums of pairs.
Thus, the arithmetic progression $A=\{x+ky, \,k=1, 2, \ldots,
N\}$ of length $N$ has the property that the image
$B=\varphi(A)$ of $A$ under $\varphi$ is a set with distinct
sums of pairs, and $\varphi$ maps bijectively $A$ onto $B$. The
construction below allows us to accumulate this effect.

For an arbitrary set $M\subseteq\mathbb R$, we define the set
$\gamma(M)$ by\footnote{Throughout the paper we use the
following notation: if $E, F \subseteq \mathbb R$ and
$a\in\mathbb R$, then $E+F=\{x+y : x\in E, \,y\in F\}$ and
$aE=\{ax : x\in E\}$.}
$$
\gamma(M)=\{0\}\cup M\cup(M+M)\cup(M-M)\cup(M+M-M).
$$
By the assumption of the theorem, $E(\varphi)$ has nonzero
measure. Let $N\geq 4$. Let $M\subseteq\mathbb R$ be a finite
set. Using Lemma~1 with $S=\gamma(M)$, we can find $x, y\in
\mathbb R$ such that
$$
\sum_{j=1}^4v_j\varphi(x+k_jy)\notin\gamma(M)
$$
for all vectors $v=(v_1, v_2, v_3, v_4)\in V$ and $k=(k_1, k_2,
k_3, k_4)\in [1, N]^4\cap\mathbb Z_0^4$. This implies that the
arithmetic progression $A=\{x+ky, \,k=1, 2, \dots, N\}$ of
length $N$ has the following properties:

(a) $\varphi$ maps bijectively $A$ onto
$B=\varphi(A)$;

(b) $B\cap M=\varnothing$;

(c) if $M$ is a set with distinct sums of pairs, then so is
$B\cup M$.

\noindent To see this we essentially repeat the above argument.
Indeed, since $0\in \gamma(M)$ and $v^3\in V$, we obtain (a).
Since $M\subseteq\gamma(M)$ and $v^7\in V$, we obtain (b). One
can easily verify assertion (c) as well; for instance, the
relation $b_1+b_2=m_2+m_3$ cannot hold for $b_1, b_2\in B$ and
$m_1, m_2\in M$, because $M+M\subseteq\gamma(M)$ and $V$
contains $v^4$ and $v^5$. The further routine verification is
left to the reader.

Using this observation, we inductively construct a family of
arithmetic progressions $A_n$, $n=4, 5, \dots$, where $A_n$ is
of length $n$, as follows. Applying Lemma~1 with $S=\{0\}$, we
find an arithmetic progression $A_4$ of length~$4$ such that
$\varphi$ maps bijectively $A_4$ onto a set $B_4$ with distinct
sums of pairs. Suppose that arithmetic progressions $A_n$,
where $A_n$ is of length $n$, are already constructed for $n=4,
5, \dots, N$. Setting $M=B_4\cup B_5\cup\dots\cup B_N$, we find
an arithmetic progression $A_{N+1}$ of length $N+1$ such that
(see (a), (b), and (c))

(i) $\varphi$ maps bijectively $A_{N+1}$ onto
$B_{N+1}$;

(ii) $B_{N+1}$ does not intersect $B_4\cup B_5\cup\dots\cup
B_N$;

(iii) $B_4\cup B_5\cup\dots\cup B_N\cup B_{N+1}$ is a set
with distinct sums of pairs.

\noindent Proceeding, we obtain $A_n$ and $B_n$ for all $n=4,
5, \dots$\,. From the construction (see (i)--(iii)) it follows
that $\varphi$ maps bijectively $\bigcup_{n=4}^\infty A_n$ onto
$\bigcup_{n=4}^\infty B_n$, and the former set contains
arbitrarily long arithmetic progressions whereas the latter is
a set with distinct sums of pairs.

\quad

\begin{center}
\textbf{4. Proof of the main Lemma}
\end{center}

To prove the main Lemma (Lemma 1) we need Lemmas 2 and 3 below.

Given a vector $k=(k_1, k_2, k_3,k_4)\in\mathbb Z^4$, a vector
$v=(v_1, v_2, v_3, v_4)\in\mathbb R^4$,  and a number
$s\in\mathbb R$, we define a set $Q(k, v, s)$ by
$$
Q(k, v, s)=\bigg\{(x, y)\in\mathbb R^2 : \sum_{j=1}^4v_j\varphi(x+k_jy)=s\bigg\}.
\eqno(3)
$$

By an interval in $\mathbb R$ we always mean a nonempty
interval.

\quad

\textsc{Lemma 2.} \emph{Suppose that the assumptions of Lemma~1
hold. Then, for any intervals $I, J\subseteq\mathbb R$, there
exist intervals $I'\subseteq I$ and $J'\subseteq J$, vectors
$k\in [1, N]^4\cap\mathbb Z_0^4$ and $v\in V$, and a number
$s\in S$ such that $I'\times J'\subseteq Q(k, v, s)$.}

\quad

One easily proves this lemma as follows. Under the assumptions
we have
$$
\bigcup_{k\in [1, N]^4\cap\mathbb Z_0^4; \,\,v\in V, \,\,s\in S} Q(k, v, s)=\mathbb R^2.
$$
So,
$$
\bigcup_{k\in [1, N]^4\cap\mathbb Z_0^4; \,\,v\in V, \,\,s\in S} Q(k, v, s)\cap(I\times J)=I\times J.
$$
Since $\varphi$ is continuous, all sets $Q(k, v, s)$ are
closed. Without loss of generality we can assume that $I$ and
$J$ are closed. Applying the Baire category theorem we complete
the proof of Lemma 2.

\quad

Given a set $F\subseteq\mathbb R$, we denote its closure by
$\overline{F}$. If $F$ is measurable, then we use $F^\circ$ to
denote the set of points of density of $F$, i.e., the set of
all $x\in F$ satisfying
$$
\frac{|F\cap I(x, \delta)|}{|I(x, \delta)|}\rightarrow 1
\qquad \textrm{as}\quad\delta\rightarrow+0.
$$
Here $|X|$ stands for the (Lebesgue) measure of a measurable
set $X\subseteq\mathbb R$ and $I(x, \delta)=(x-\delta,
x+\delta)$. As is known, almost all points of a measurable set
are its points of density.

\quad

\textsc{Lemma 3.} \emph{Let $E_1$ and $E_2$ be measurable sets
in $\mathbb R$ such that
$\overline{E_1^\circ}\cap\overline{E_2^\circ}$ has nonzero
Lebesgue measure. Then, for any positive integer $N$ and any
partition of the set $\{1, 2, \dots, N\}$ into two disjoint
sets $\mathcal{E}_1$ and $\mathcal{E}_2$, there exists an
arithmetic progression $t_k=x+ky$, $k=1, 2, \dots, N$, of
length $N$ such that $t_k\in E_1^\circ$ for $k\in
\mathcal{E}_1$ and $t_k\in E_2^\circ$ for $k\in
\mathcal{E}_2$.}

\quad

This combinatorial lemma was obtained by the author and A.
Olevskii in [9, Lemma~1] (see also [8]). It plays one of the
key roles in the proof of the $M_p$ -version of the
Beurling--Helson theorem (see Sec. 5, Remark 5). In a slightly
weaker form it was obtained earlier by the same authors in the
work [7] on idempotent Fourier multipliers.

\quad

We now proceed directly to the proof of Lemma~1 (the main
Lemma). The proof is split into two steps.

\emph{Step 1.} First, we  obtain a weaker result;
namely, we show that, under the assumptions of Lemma~1, the set
$E(\varphi)$ is nowhere dense.

Let $\Delta=(a, b)$ be an arbitrary interval in $\mathbb R$.
Let us show that it contains a subinterval on which $\varphi$
is linear; this will prove our claim. To this end consider the
following two intervals $I$ and $J$:
$$
I=\bigg(a, \,\frac{a+b}{2}\bigg), \qquad J=\bigg(0, \,\frac{b-a}{2N}\bigg).
\eqno(4)
$$
Using Lemma~2, we find vectors $v=(v_1, v_2, v_3, v_4)\in
V$ and $k=(k_1, k_2, k_3, k_4)\in [1, N]^4\cap\mathbb Z_0^4$,
a number $s\in S$, and intervals $I'\subseteq I$ and $J'\subseteq J$
such that $I'\times J'\subseteq Q(k, v, s)$, i.e.,
$$
\sum_{j=1}^4v_j\varphi(x+k_jy)=s \qquad\textrm{for all} \quad (x, y)\in I'\times J'.
\eqno(5)
$$

Choose an infinitely differentiable nonnegative function $\rho$
on $\mathbb R$ so that $\mathrm{supp}\,\rho\subseteq[-1, 1]$
and $\int_\mathbb R \rho(x)dx=1$. For each $\varepsilon>0$, we
set
$\rho_\varepsilon(t)=\frac{1}{\varepsilon}\rho(\frac{t}{\varepsilon})$.
We have $\mathrm{supp}\,\rho_\varepsilon\subseteq[-\varepsilon,
\varepsilon]$. Consider the convolution
$\varphi_\varepsilon=\varphi\ast\rho_\varepsilon$:
$$
\varphi_\varepsilon(u)=\int_\mathbb R \varphi(u-t)\rho_\varepsilon(t)dt.
$$
Obviously $\varphi_\varepsilon$ is infinitely differentiable.
Note also that since $\varphi$ is continuous, it follows that
$\varphi_\varepsilon$ converges pointwise to $\varphi$ as
$\varepsilon\rightarrow +0$.

Let $I''$ be the interval concentric with $I'$ and of length
three times smaller than that of $I'$. We set $\varepsilon_0$
equal to the length of $I''$. For any $t$ with
$|t|<\varepsilon_0$ and any point $(x, y)\in I''\times J'$, we
have  $(x-t, y)\in I'\times J'$; therefore (see (5)),
$$
\sum_{j=1}^4v_j\varphi(x-t+k_jy)=s.
$$
Thus, we see that
$$
\sum_{j=1}^4v_j\varphi_\varepsilon(x+k_jy)=s \qquad\textrm{for all}
\quad (x, y)\in I''\times J' \quad \textrm{and} \quad 0<\varepsilon<\varepsilon_0.
$$

Differentiating this relation three times, namely, taking the
derivatives $\frac{\partial^3}{\partial x^{3-q}\partial
y^{q}}$, $0\leq q\leq 3$, we obtain
$$
\sum_{j=1}^4 k_j^q v_j\varphi_\varepsilon^{'''}(x+k_jy)=0,
\qquad 0\leq q\leq 3, \qquad x\in I'', \,y\in J'.
$$
Since $k_j$, $1\leq j\leq 4$, are pairwise different (recall
that $k\in\mathbb Z_0^4$), it follows that the matrix
$\{k_j^q\}_{1\leq j\leq 4, \,0\leq q\leq 3}$ has nonzero
determinant. Since not all $v_j$ vanish, we see that there
exists a $j_0$ (we can take, e.g., $j_0=1$, see (2)) such that
$\varphi_\varepsilon^{'''}(x+k_{j_0}y)=0$ for all $x\in I''$
and $y\in J'$. Thus, if $0<\varepsilon<\varepsilon_0$, then
$\varphi_\varepsilon^{'''}(t)=0$ for all $t\in I''+k_{j_0}J'$.
Hence, $\varphi_\varepsilon$ is a polynomial of degree at
most~2 on the interval $\widetilde{\Delta}=I''+k_{j_0}J'$.
Letting $\varepsilon\rightarrow+0$, we see that $\varphi$
coincides with a polynomial $P$ of degree at most $2$ on
$\widetilde{\Delta}$. Since $k\in [1, N]^4$, we have $1\leq
k_{j_0}\leq N$, whence
$\widetilde{\Delta}=I''+k_{j_0}J'\subseteq I+k_{j_0}J\subseteq
\Delta$ (see(4)).

Let us show that the degree of the polynomial $P$ is strictly
less than~$2$. We repeat part of the argument used above, this
time for the interval $\widetilde{\Delta}=~(\widetilde{a},
\widetilde{b})$ instead of $\Delta=(a, b)$. Namely, we consider
the following intervals $\widetilde{I}$ and $\widetilde{J}$:
$$
\widetilde{I}=\bigg(\widetilde{a}, \,\frac{\widetilde{a}+\widetilde{b}}{2}\bigg),
\qquad \widetilde{J}=\bigg(0, \,\frac{\widetilde{b}-\widetilde{a}}{2N}\bigg).
$$
Using Lemma~2, we find vectors $v=(v_1, v_2, v_3, v_4)\in
V$, $k=(k_1, k_2, k_3, k_4)\in [1, N]^4\cap\mathbb Z_0^4$,
a number $s\in S$, and intervals $\widetilde{I}'\subseteq
\widetilde{I}$ and $\widetilde{J}'\subseteq \widetilde{J}$ such
that $\widetilde{I}'\times\widetilde{J}'\subseteq Q(k, v, s)$,
i.e.,
$$
\sum_{j=0}^4v_j\varphi(x+k_j y)=s \qquad\textrm{for all}
\quad (x, y)\in\widetilde{I}'\times\widetilde{J}'.
\eqno(6)
$$
Since all points $x+k_jy$ with $x\in\widetilde{I}'$, $y\in
\widetilde{J}'$, and $j=1, 2, 3, 4,$ are in
$\widetilde{\Delta}$ and since $\varphi$ coincides with $P$ on
$\widetilde{\Delta}$, from (6) it follows that
$$
\sum_{j=0}^4v_j P(x+k_j y)=s \qquad\textrm{for all}
\quad (x, y)\in\widetilde{I}'\times\widetilde{J}'.
\eqno(7)
$$
Assuming that the degree of $P$ equals $2$, we have
$P''\equiv\mathrm{const}\neq 0$. Twice differentiating relation
(7), that is, taking the derivatives
$\frac{\partial^2}{\partial^{2-q} x\partial^q y}$, $0\leq q\leq
2$, we see that
$$
\sum_{j=0}^4 k_j^q v_j=0, \qquad 0\leq q\leq 2.
$$
Thus, the vector $v$ belongs to the kernel of the matrix
$$
M(k)=\left(
       \begin{array}{cccc}
         1 & 1 & 1 & 1 \\
         k_1 & k_2 & k_3 & k_4 \\
         k_1^2 & k_2^2 & k_3^2 & k_4^2 \\
       \end{array}
     \right).
$$
It remains to observe that this is impossible: none of the
vectors of the family $V$ (see (2)) belongs to $\ker M(k)$
whenever $k_1$, $k_2$, $k_3$, $k_4$ are pairwise different
positive integers. The verification is left to the reader. Thus
we see that $P$ is a linear function. This completes the proof
of our claim that $E(\varphi)$ is nowhere dense.

\emph{Step 2.} Now, we show that, under the assumptions of
Lemma 1, the set $E(\varphi)$ has measure zero. The core of
this step is the combinatorial Lemma~3.

Let $\Omega$ denote the family of all open intervals
complementary to $E(\varphi)$, i.e., the family of connected
components of the complement $\mathbb R\setminus E(\varphi)$
(recall that the set $E(\varphi)$ is closed). For every
interval $I\in\Omega$, we have $\varphi(t)=P_I(t), \,t\in I,$
where $P_I$ is a linear function.

Now, suppose that, contrary to the assertion of Lemma~1, the
set $E(\varphi)$ has nonzero measure. Let $E$ be the set of
accumulation points of $E(\varphi)$. Under assumption that
$E(\varphi)$ has nonzero measure, the same holds for $E$.

Note that if $x_0\in E$, then (since $E(\varphi)$ is closed and
nowhere dense) any neighborhood of $x_0$ contains infinitely
many intervals $I\in\Omega$ with the property that the
corresponding functions $P_I$ are pairwise different. Indeed,
otherwise, the point $(x_0, \varphi(x_0))\in\mathbb R^2$ has a
neighborhood $J$ such that the piece $J\cap G$ of the graph $G$
of $\varphi$ is contained in a finite union of straight lines,
which is possible only if $x_0$ is an isolated point of
$E(\varphi)$ or does not belong to $E(\varphi)$ at all.

For each $n=1, 2, \dots$, consider the open $1/n$ -neighborhood
of the set $E$. Let $\Delta_k^n$, $k=1, 2, \dots,$ be the
family of all connected components of this neighborhood. We
renumber the intervals $\Delta_k^n$, $n=1, 2, \dots$, $k=1, 2,
\dots$, as $\Delta_1, \Delta_2, \dots$\,. Each of the intervals
$\Delta_j$ contains an accumulation point of the set
$E(\varphi)$. Hence, it contains infinitely many intervals $I$
complementary to $E(\varphi)$ with the property that the
corresponding functions $P_I$ are pairwise different. We choose
an interval $I_1\in\Omega$ contained in $\Delta_1$. Having
chosen intervals $I_1, I_2, \dots, I_j\in\Omega$ contained in
$\Delta_1, \Delta_2, \dots, \Delta_j$, respectively, we choose
an interval $I_{j+1}\in\Omega$ so that
$I_{j+1}\subseteq\Delta_{j+1}$ and none of nontrivial linear
combinations of $P_{I_{j+1}}, P_{I_1}, P_{I_2}, \dots, P_{I_j}$
with coefficients $0, \pm 1, \pm 2$ can be identically equal to
$s$ whenever $s\in S$ (by a nontrivial linear combination we
mean a combination not all of whose coefficients are zero).
Clearly, such an interval always exists, because $S$ is finite
and there are only finitely many linear combinations of the
functions $P_{I_1}, P_{I_2}, \dots, P_{I_j}$ with coefficients
$0, \pm 1, \pm 2, \pm1/2$. Proceeding by induction, we obtain
intervals $I_m\in\Omega$, $m=1, 2, \dots,$ with the following
two properties: firstly, the intervals $I_m$ accumulate to $E$,
i.e., any neighborhood of any point of $E$ contains an interval
that belongs to the family $\{I_m\}$, and, secondly, no
nontrivial linear combination of the corresponding linear
functions $P_{I_m}$, $m=1, 2, \dots$, with coefficients $0, \pm
1, \pm 2$ can be identically equal to $s$ whenever $s\in S$.

Denote the union of the intervals $I_m$, $m=1, 2, \dots$,
by~$U$. Observe that the sets $E^\circ$ and $U^\circ$ (the sets
of points of density of $E$ and $U$, respectively) satisfy
$\overline{E^\circ}\cap\overline{U^\circ}\supseteq E^\circ$.
Hence, $\overline{E^\circ}\cap\overline{U^\circ}$ has nonzero
measure. Using Lemma~3, we find an arithmetic progression
$t_k=a+kd$, $k=1, 2, \dots, 2N$, of length $2N$ such that its
terms with odd indices belong to $E$ and those with even
indices belong to $U$. Consider only the terms with even
indices. Clearly they form a progression $x_0+ky_0$, $k=1, 2,
\dots, N$, of length $N$ with the property that no two
different terms of this progression belong to the same interval
of the family $\{I_m\}$. We consider now only those intervals
of the family $\{I_m\}$ which contain a point of this
progression. For $k=1, 2, \ldots, N$ denote the interval of the
family $\{I_m\}$ that contains the point $x_0+ky_0$ by $H_k$.
Let $P_k$, $k=1, 2, \dots, N$, be the corresponding linear
functions, i.e., $P_k=P_{H_k}$.

Clearly, if $I$ and $J$ are sufficiently small neighborhoods of
the points $x_0$ and $y_0$, respectively, then for all $x\in I$
and $y\in J$ we have $x+ky\in H_k, \,k=1, 2, \dots, N$. We fix
these $I$ and~$J$.

Again applying Lemma~2, we find vectors $v=(v_1, v_2, v_3,
v_4)\in V$ and $k=(k_1, k_2, k_3, k_4)\in[1, N]^4\cap\mathbb
Z_0^4$, a number $s\in S$, and intervals $I'\subseteq I$ and
$J'\subseteq J$ such that
$$
\sum_{j=1}^4 v_j \varphi(x+k_j y)=s \quad\textrm{for all}\quad (x, y)\in I'\times J'.
$$
This implies
$$
\sum_{j=1}^4 v_j P_{k_j}(x+k_j y)=s \quad\textrm{for all}\quad (x, y)\in I'\times J'.
$$
Clearly, if an affine function of two variables identically
equals $s$ on a rectangle in $\mathbb R^2$, then it identically
equals $s$ in the entire plane~$\mathbb R^2$. Thus,
$$
\sum_{j=1}^4 v_j P_{k_j}(x+k_j y)=s, \quad (x, y)\in\mathbb R^2.
$$
Setting $y=0$, we see that
$$
\sum_{j=1}^4 v_j P_{k_j}(x)=s, \qquad x\in\mathbb R,
$$
which is impossible because the coordinates of each vector
$v\in V$ are $0, \pm 1$ or $\pm 2$ and not all of them are zero
(see(2)). This proves Lemma~1 and, thereby, the theorem.

\quad

\begin{center}
\textbf{5. Remarks}
\end{center}

1. The theorem proved in this paper admits a generalisation for
mappings $\varphi \colon I\rightarrow\mathbb R$, where $I$ is
an interval in $\mathbb R$. Indeed, without loss of generality,
we can assume that $I$ is a closed interval, $I=[a, b]$. It
suffices to consider a continuous extension of $\varphi$
constant on the rays $(-\infty, a)$ and $(b, +\infty)$, and
apply the original version of the theorem to the extension.

2. The following assertion on affine copies of $\mathbb Z$
holds: \emph{If $\varphi \colon \mathbb R\rightarrow\mathbb R$
is a continuous nowhere linear (that is, nonlinear on every
interval) mapping, then there exists an affine copy $a\mathbb
Z+b$ of $\mathbb Z$ such that $\varphi$ maps bijectively
$a\mathbb Z+b$ onto a certain set with distinct sums of pairs.}
This can be proved by a modification of the first step in the
proof of the main Lemma (without use of the second step based
on Lemma~3). Indeed, let $V'$ be a family of the first three
vectors $v^1, v^2, v^3$ defined in (2). Assuming that the
assertion on copies is not true, we have (see (3))
$$
\bigcup_{k\in\mathbb Z_0^4; \,\,v\in V'} Q(k, v, 0)=\mathbb R^2.
$$
Using categorical considerations (similar to those in the proof
of Lemma~2) we obtain two intervals $I, J\subseteq\mathbb R$, a
vector $k=(k_1, k_2, k_3, k_4)\in\mathbb Z_0^4$, and a vector
$v=(v_1, v_2, v_3, v_4)\in V'$ such that
$$
\sum_{j=1}^4v_j\varphi(x+k_jy)=0 \qquad\textrm{for all} \quad (x, y)\in I\times J.
$$
Choosing an interval $I'$ concentric with $I$ and of strictly
smaller length than that of $I$, and repeating the argument of
the first step in the proof of the main Lemma, we obtain that
$\varphi$ is linear on $I'+k_{j_0}J$, which contradicts the
assumption on~$\varphi$.

3. Consider the algebra $B(\mathbb R)$ of Fourier transforms of
measures on $\mathbb R$. According to the well-known
Beurling--Helson theorem [2], if $\varphi$ is a real continuous
function on $\mathbb R$ such that
$$
\|e^{in\varphi}\|_{B(\mathbb R)}=O(1), \quad n\in\mathbb Z,
\eqno(8)
$$
then $\varphi$ is linear. This theorem has a version for the
Wiener algebra $A(\mathbb T)$ of absolutely convergent Fourier
series on the circle $\mathbb T$, which is due to Kahane [5,
Ch.~VI]. In particular, the Beurling--Helson theorem implies
that only trivial (i.e., linear) changes of variable are
allowable in $B(\mathbb R)$ (the same holds for $A(\mathbb
T)$). An essential point in the proof of the Beurling--Helson
theorem is the observation that condition (8) implies that
$\varphi$ cannot map bijectively long arithmetic progressions
onto sets which are $\mathrm{mod}~2\pi$~-independent over
integers. (Subsequently, Kahane [6] gave a proof based on
completely different argument.) So, the question of how thin
the images of arithmetic progressions under continuous mappings
can be traces back to Beurling and Helson.

4. A closed set $E\subseteq\mathbb R$ is called a Helson set if
any continuous function on $E$ vanishing at infinity is the
restriction to $E$ of the Fourier transform of a function in
$L^1(\mathbb R)$. Equivalently (see [4]), $E$ is a Helson set
if, given a measure $\mu$ on $\mathbb R$, we have
$\|\mu\|_{M(\mathbb R)}\leq
c\|\widehat{\mu}\|_{L^\infty(\mathbb R)}$, where
$\|\mu\|_{M(\mathbb R)}$ is the variation of $\mu$ and $c>0$
does not depend on $\mu$. In particular, taking for $\mu$ a
linear combination of point masses, we see that
$$
\sum_{k=1}^N |c_k|\leq c
\bigg\|\sum_{k=1}^N c_k e^{i\lambda_kt}\bigg\|_{L^\infty(\mathbb R)}
$$
for any $N$, any pairwise different $\lambda_1, \lambda_2,
\dots, \lambda_N\in E$, and any (complex) numbers $c_k$, $k=1,
2, \dots, N$. Countable Helson sets are called Sidon sets. In
[3], Graham considered self-mappings of $\mathbb R$ which take
Sidon sets to Sidon sets. He called a mapping $\varphi$
countably piecewise affine if the set of $x$ such that
$\varphi$ is affine in a neighborhood of $x$ is dense in
$\mathbb R$ (which in our terms is the same as to say that
$E(\varphi)$ is nowhere dense). Graham showed that if a
self-homeomorphism $h$ of $\mathbb R$ has the property that the
image $h(E)$ of any Sidon set $E$ is a Sidon set, then $h$ is
countably piecewise affine with a finite number of slopes. The
use of categorical argument and subsequent use of convolution
in the proof of the Theorem of the present paper were suggested
by this work of Graham. We also note that Graham conjectured
that his result on Sidon sets can be supplemented by the
assertion that the set of nonlinearity of $h$ has measure zero.
It is very plausible that this is indeed the case and perhaps
it can be proved by an argument similar to that used at the
second step of the proof of Lemma 1 of the present paper.

5. Mappings $\varphi$ whose sets of nonlinearity have measure
zero also appeared in relation with analogues of the
Beurling--Helson theorem for the algebras $M_p(\mathbb R)$ of
Fourier multipliers (it is well-known that $M_1(\mathbb
R)=M_\infty(\mathbb R)=B(\mathbb R)$ and the corresponding
norms coincide). As it turned out [9] (see also [8]), if
$\varphi \colon \mathbb R\rightarrow [0, 2\pi[$ is a measurable
function such that $\|e^{in\varphi}\|_{M_p(\mathbb R)}=O(1),
\,n\in\mathbb Z,$ for some $p$, $1<p<\infty$, $p\neq 2$, then
$E(\varphi)$ has measure zero and the set of distinct slopes of
$\varphi$ is finite.

\begin{center}
\textbf{References}
\end{center}

\flushleft
\begin{enumerate}

\item A. S. Besicovitch, \emph{Almost Periodic Functions},
    Cambridge University Press, Cambridge, 1932.

\item  A. Beurling, H. Helson, ``Fourier--Stieltjes
    transforms with bounded powers,''  \emph{Math. Scand.},
    \textbf{1} (1953), 120--126.

\item C. C. Graham, ``Mappings that preserve Sidon sets in
    $R$,'' \emph{Arkiv f\"or Matematik}, \textbf{19}:1
    (1981), 217--221.

\item  C. C. Graham, O. C. McGehee, \emph{Essays in
    Commutative Harmonic Analysis,} Grundlehren der
    Mathematischen Wissenschaften, Vol. 238, Springer-Verlag,
    New York, 1979.

\item J.-P. Kahane, \emph{S\'erie de Fourier absolument
    convergentes}, Springer-Verlag, Berlin--Heidelberg--New
    York, 1970.

\item J.-P. Kahane, ``Quatre le\c cons sur les
    hom\'eomorphismes du circle et les s\'eries de Fourier,''
    in: \emph{Topics in Modern Harmonic Analysis, Vol. II,
    Ist. Naz. Alta Mat. Francesco Severi, Roma}, 1983,
    955--990.

\item V. Lebedev and  A. Olevski\v{\i}, ``Idempotents of
    Fourier multiplier algebra,'' \emph{Geometric and
    Functional Analysis (GAFA)}, \textbf{4}:5 (1994),
    539--544.

\item V. Lebedev and  A. Olevski\v{\i}, ``Bounded groups of
    translation invariant operators,'' \emph{C. R. Acad. Sci.
    Paris, S\'er. \emph{I} Math.}, \textbf{322}:2 (1996),
    143--147.

\item V. V. Lebedev, A. M. Olevskii, ``$L^p$ -Fourier
    multipliers with bounded powers,'' \emph{Izvestya:
    Mathematics}, \textbf{70}:3 (2006), 549--585.

\item B. M. Levitan, \emph{Almost Periodic Functions},
    Gostekhizdat, Moscow, 1953 (in Russian).

\item W. Rudin, ``Trigonometric series with gaps,''
    \emph{Journal of Mathematics and Mechanics}, \textbf{9}:2
    (1960), 203--227.

\end{enumerate}

\quad

\quad

\qquad School of Applied Mathematics\\
\qquad National Research University Higher School of Economics\\
\qquad 34, Tallinskaya St.\\
\qquad Moscow, 123458 Russia\\
\qquad E-mail address: \emph{lebedevhome@gmail.com}

\end{document}